\documentclass{article}

\usepackage{amssymb}
\usepackage{mathrsfs}
\usepackage{epsfig}

\begin{document}

\newtheorem{thm}{Theorem}
\newtheorem{definition}[thm]{Definition}
\newtheorem{lem}[thm]{Lemma}
\newtheorem{prop}[thm]{Proposition}
\newtheorem{cor}[thm]{Corollary}
\newtheorem{rmk}[thm]{Remark}

\newcommand{\qed}{\hfill \mbox{\raggedright \rule{.1in}{.1in}}}

\title{Strong accessibility for hyperbolic groups}
\author{Diane M. Vavrichek\footnote{Research supported in part by NSF grant DMS-0602191.}}
\date{}
\maketitle

\begin{abstract}
We use an accessibility result of Delzant and Potyagailo to prove Swarup's Strong Accessibility Conjecture for Gromov hyperbolic groups with no 2-torsion.  
\end{abstract}


A particular notion of a hierarchy for 3-manifolds motivated Swarup to make the following conjecture.  
\\

\noindent {\bf Swarup's Strong Accessibility Conjecture} {\it
Let $G$ be a hyperbolic group.  Decompose $G$ maximally over finite subgroups, and 
then take the resulting vertex groups, and decompose those maximally over two-ended 
subgroups.  Now repeat this process on the new vertex groups and so on.  Then this 
process must eventually terminate, with subgroups of $G$ that are unsplittable over 
finite and two-ended subgroups.
} \\

We prove this result in the case that $G$ has no 2-torsion, making heavy use of the main result from \cite{DP}.
In \cite{DP}, Delzant and Potyagailo show the existence of finite hierarchies for finitely presented groups with no 2-torsion, over so-called elementary families of subgroups.
Even ignoring the cases when $G$ has 2-torsion, the conjecture of Swarup is not a special case of this accessibility result for two reasons.  Firstly, to get a finite hierarchy from the construction in \cite{DP}, we are required to choose one elementary family of subgroups to decompose over, and stick with it for the entire process.  In contrast, for Swarup's conjecture, we must alternate between decomposing over finite and two-ended subgroups.
Secondly, given a group and a family of elementary subgroups, \cite{DP} shows the existence of one finite hierarchy over the family.  For Swarup's conjecture, it must instead be shown that any hierarchy as described is finite.
By analyzing equivariant maps between $G$-trees, we are able to overcome these difficulties to prove the conjecture in the case in which $G$ has no 2-torsion.
As a corollary to this, we get the following result about finite hierarchies in 3-manifolds:
\\

\noindent {\bf Theorem \ref{3mans}} {\it
Let $M$ be an irreducible, orientable, compact 3-manifold with hyperbolic fundamental group.
The process of decomposing $M$ along any maximal, disjoint collection of compressing 
disks, then decomposing the resulting manifolds along maximal, disjoint collections of 
essential annuli, then the resulting manifolds along compressing disks, then again along 
essential annuli and so on, must eventually terminate with a collection of manifolds 
which have incompressible boundary and admit no essential annuli, or are 3-balls.
}\\


To begin, we now review some basic notions and terminology that we will need, mainly coming from Bass-Serre theory.
Let a group $G$ act simplicially on the left on a simplicial tree $\tau$, and let the action be without inversions, i.e., such that no element of $G$ fixes an edge of $\tau$, but swaps its vertices.  Then we say that $\tau$ is a {\it $G$-tree}.  

Label each vertex $v_0$ of $\Gamma = G\backslash \tau$ with the stabilizer $V$ of one of its preimages $v$ under the projection map $\tau \to G\backslash \tau$, and label each edge $e_0$ with the stabilizer $E$ of one of its preimages as well.
Note that these stabilizers, or {\it vertex} and {\it edge groups}, are uniquely determined up to conjugacy.
Then, to each pair $(e_0, v_0)$ of a vertex and incident edge in $\Gamma$, associate an injective homomorphism from $E$ to $V$ induced by the inclusion of the stabilizer of a lift of $e_0$ into the stabilizer of a lift of $v_0$.
Call $\Gamma$, together with all this data, a {\it graph of groups} structure for $G$, and denote the graph and data also by $\Gamma$.
We shall also say that $\Gamma$ is a {\it decomposition} of $G$.

$G$ is said to act {\it minimally} on a $G$-tree $\tau$ if $\tau$ contains no proper, $G$-invariant subtrees.
If $G$ acts minimally on $\tau$, $\tau$ is not one vertex, and $\Gamma$ is finite, then we call $\Gamma$ a {\it proper decomposition} of $G$.
Note that if $G$ is finitely generated, acts minimally on $\tau$, and $\tau$ is not a vertex, then $\Gamma = G \backslash \tau$ must be a finite graph, and hence is a proper decomposition  of $G$.

Any simplicial, $G$-equivariant map between $G$-trees induces a {\it graph of groups map} of the associated decompositions of $G$, i.e. a simplicial, surjective map of the underlying graphs, that may collapse edges to vertices, and injective homomorphisms of the vertex and edge groups which commute with the injections from edge groups into their incident vertex groups.  Let a $G${\it -map} be a simplicial, $G$-equivariant map between two $G$-trees, that does not collapse any edge to a vertex.  

We will now take $G$ to be finitely generated.  Let $\tau$ and $\tau'$ be $G$-trees which are not vertices and on which $G$ acts minimally, and let $\Gamma = G \backslash \tau$ and $\Gamma' = G \backslash \tau'$.  
Assume that there is a graph of groups map $\Gamma' \to \Gamma$, that is induced from a map $\tau' \to \tau$ that is not a simplicial homeomorphism.  Then
we call the decomposition $\Gamma'$ a {\it refinement} of $\Gamma$.  If, for each edge $e$ with vertices $x$ and $y$ of $\tau'$ which is collapsed to a vertex of $\tau$, either $x$ and $y$ are in the same $G$-orbit, or $X=$ stab$(x)$ and $Y=$ stab$(y)$ properly contain $E = $ stab$(e)$, then we call $\Gamma'$ a {\it proper refinement} of $\Gamma$.

If all edge stabilizers of a decomposition $\Gamma$ of $G$ are in some family $\mathscr{C}$ of subgroups of $G$, then we say that $\Gamma$ is a decomposition of $G$ {\it over} $\mathscr{C}$.  
Since the edge groups of $\Gamma$ are determined only up to conjugacy, $\mathscr{C}$ should be closed under conjugacy.
Note that if $\Gamma$ is a decomposition of $G$ over $\mathscr{C}$, and $\Gamma'$ is a refinement of $\Gamma$, then $\Gamma'$ is a decomposition of $G$ over the elements of $\mathscr{C}$ and their subgroups.

A decomposition of $G$ with one edge is a {\it splitting} of $G$, and a proper decomposition of $G$ with one edge is a {\it proper splitting} of $G$.  If there exist no proper splittings of $G$ over a family $\mathscr{C}$ as above, then we say that $G$ is {\it unsplittable over} $\mathscr{C}$.  

Note that if $G$ admits a proper decomposition $\Gamma'$ over $\mathscr{C}$, arising from an action on a $G$-tree $\tau'$, then for any edge $e$ of $\Gamma'$ with edge group $E$, there is a proper splitting $\Gamma$ of $G$ associated to $e$, where $\Gamma$ has one edge with edge group $E$, and $\Gamma'$ is a refinement of $\Gamma$.  To see this, let $e$ be an edge in $\tau'$ with stabilizer $E$.  Then, if $G\cdot $int$(e)$ denotes the $G$-orbit of the interior of $e$ in $\tau'$,  let $\tau$ be the $G$-tree obtained by collapsing the components of $\tau' -G\cdot $int$(e)$ to vertices, with the action of $G$ induced from the action of $G$ on $\tau'$.  Then we may take $\Gamma$ to be $G \backslash \tau$.

If $G$ has a decomposition $\Gamma$, and the vertex group of a vertex $v$ of $\Gamma$ admits a splitting, then we say that the splitting is {\it compatible} with the decomposition if there exists a refinement of $\Gamma$ in which $v$ is replaced with an edge that is associated to the splitting, as in the previous paragraph.  Equivalently, the splitting is compatible with the decomposition if a conjugate of each edge group of the edges incident to $v$ is contained in a vertex group of the splitting.

Consider a group $G$, and a family $\mathscr{C}$ of subgroups of $G$ which is closed under conjugacy.  A {\it hierarchy} for $G$ over $\mathscr{C}$ is a sequence $\mathscr{G}_0, \mathscr{G}_1, \mathscr{G}_2, \ldots$ of finite sets of conjugacy classes of subgroups of $G$, defined inductively as follows.
The set $\mathscr{G}_0$ contains only (the conjugacy class of) $G$.  If $i>0$, then for any conjugacy class in $\mathscr{G}_{i-1}$, either $\mathscr{G}_i$ contains that conjugacy class, or $\mathscr{G}_i$ contains the conjugacy classes of the vertex groups of some proper decomposition of a representative of that class over $\mathscr{C}$.  
We require that at least one representative from $\mathscr{G}_{i-1}$ be decomposed. 
(Since the structures of decompositions of representatives of conjugacy classes over $\mathscr{C}$ are invariant under conjugacy, our choices of representatives do not matter.)  

If this process terminates with some $\mathscr{G}_N$ which contains only conjugacy classes of subgroups which are unsplittable over $\mathscr{C}$, then we say that the hierarchy is {\it finite}.  

We note that the existence of a finite hierarchy over $\mathscr{C}$ does not, in general, imply the existence of any kind of maximal decomposition over $\mathscr{C}$, since the splittings of vertex groups need not be compatible with the decompositions producing those vertex groups.

If $\Gamma$ is a decomposition of $G$ over $\mathscr{C}$, and if there exists no proper refinement of $\Gamma$ over $\mathscr{C}$, then we say that $\Gamma$ is a {\it maximal decomposition} of $G$ over $\mathscr{C}$.  

One could alternatively define a maximal decomposition of $G$ over $\mathscr{C}$ to be a maximal collection of compatible splittings over $\mathscr{C}$.  This is a stronger requirement than what we have made, and such collections in fact need not be finite.
For example, consider any group $G$ that has an infinite descending chain of subgroups $G \supset C_0 \supset C_1 \supset C_2 \supset \ldots$.  Then we have 
$$G = G*_{C_0}C_0 = G*_{C_0}C_0 *_{C_1}C_1 =  G*_{C_0}C_0 *_{C_1}C_1*_{C_2}C_2 = \ldots$$
Less trivially, consider the Baumslag-Solitar group $H = $ BS$(1,2) = \langle x,t: t^{-1}xt=x^2\rangle$.  The normal closure of $\langle x\rangle$ in $H$ is isomorphic to $\mathbb{Z}[\frac{1}{2}]$ under addition, by an isomorphism which takes $x$ to $1$, and  $t^ixt^{-i}$ to $\frac{1}{2^i}$ for each $i$.  Let $A_i$ denote the infinite cyclic subgroup generated by $t^{i}xt^{-i}$, and let $K = H*_{A_0}H$.  Note that $A_0\subset A_1\subset A_2\subset \ldots$, and that $K$ is finitely presented.
We can refine the given decomposition of $K$ as many times as we please, for we have that 
$$K = H*_{A_0} (A_1*_{A_1} H)  =  H*_{A_0}( A_1*_{A_1} (A_2*_{A_2} H)) =\ldots$$ 
Hence both of the examples above have sequences of 
refinements that do not terminate.


Two subgroups $H$ and $H'$ of a group $G$ are said to be {\it commensurable} if their intersection is of finite index in both.  The {\it commensurizer}, Comm$_G(H)$, of $H$ in $G$ is the subgroup of elements $g$ of $G$ such that $H$ and $gHg^{-1}$ are commensurable.


Let $G$ be a finitely generated group, with a finite generating set $S$ such that $s \in S$ implies $s^{-1} \notin S$.  The {\it Cayley graph} $\Gamma_G(S)$ of $G$ with respect to $S$ is a graph with vertex set equal to $G$, and one edge for each pair of $g \in G$ and $s \in S$, connecting the vertex $g$ with the vertex $gs$. 
By giving each edge in $\Gamma_G(S)$ length one, we may view the graph as a metric space.  We remark that, if $S$ and $S'$ are two different finite generating sets for $G$, then $\Gamma_G(S)$ and $\Gamma_G(S')$ are quasi-isometric.  

The number of {\it ends} of a locally finite simplicial complex $X$, denoted $e(X)$, is defined to be the supremum over all finite subcomplexes $K$ of $X$ of the number of infinite components of $X-K$.
Then the number of {\it ends} of a finitely generated group $G$, $e(G)$, is the number of ends of $\Gamma_G(S)$, where $S$ is a finite generating set for $G$.  This does  not depend on the choice of $S$.

A finitely generated group $G$ is said to be {\it hyperbolic} if there is some $\delta >0$ and some finite generating set $S$ for $G$ such that, for any geodesic triangle in $\Gamma_G(S)$, each side of the triangle is contained in the union of the $\delta$-neighborhoods of the other two sides. While the value of $\delta$ depends on our choice of a generating set, we note that hyperbolicity does not.


In \cite{DP}, the authors prove the existence of a finite hierarchy for any finitely presented group with no 2-torsion over any family of ``elementary'' subgroups, which are defined as follows.

\begin{definition}\label{def_of_elem}
If $G$ is a finitely presented group, and $\mathscr{C}$ a family of subgroups of $G$, 
then $\mathscr{C}$ is said to be {\rm elementary} if the following conditions are 
satisfied:
\begin{enumerate}
\item{If $C \in \mathscr{C}$, then all subgroups and conjugates of $C$ are in 
$\mathscr{C}$.}
\item{Each infinite element of $\mathscr{C}$ is contained in a unique maximal subgroup in $\mathscr{C}$.}
\item{Ascending unions of finite subgroups in $\mathscr{C}$ are contained in 
$\mathscr{C}$.}
\item{If any $C \in \mathscr{C}$ acts on a tree, then $C$ fixes a point in the tree, 
fixes a point in the boundary at infinity of the tree, or preserves but permutes two 
points in the boundary at infinity.}
\item{If $C \in \mathscr{C}$ is an  infinite, maximal element of $\mathscr{C}$ and 
$gCg^{-1} = C$, then $g \in C$.}
\end{enumerate}
\end{definition}

We will be interested in applying the results of \cite{DP}   
when $G$ is hyperbolic and $\mathscr{C}$ is the set of all finite and two-ended subgroups of $G$, so we must show the following.

\begin{prop} \label{elem}
If $G$ is a hyperbolic group, and $\mathscr{C}$ is the set of all finite and two-ended subgroups of $G$, then $\mathscr{C}$ is elementary.
\end{prop}

We now lay out several lemmas that we will use to show this.
First, we have the following classical theorem which is shown in \cite{ScW} and which characterizes two-ended groups.
\\

\begin{thm} \label{char_two-ended}
The following conditions on a finitely generated group $H$ are equivalent:
\begin{enumerate}
  \item{$e(H) = 2$,}
  \item{$H$ has an infinite cyclic subgroup of finite index,}
  \item{$H$ has a finite normal subgroup with quotient $\mathbb{Z}$ or $\mathbb{Z}_2*\mathbb{Z}_2$,}
  \item{$H = F*_F$ with $F$ finite, or $H = A*_FB$ with $F$ finite and $[A:F]=[B:F]=2$.}
\end{enumerate}
\end{thm}

Note that if $H$ is two-ended, and has a finite subgroup $F$ such that $H/F \cong \mathbb{Z}$, then $F$ must be the set of elements of $H$ that have finite order.  If $H/F \cong \mathbb{Z}_2 * \mathbb{Z}_2$, then $F$ is also characteristic:

\begin{lem} \label{F_unique}
If $H$ is a two-ended group that contains a finite normal subgroup $F$ such that 
$H/F \cong \mathbb{Z}$ or $\mathbb{Z}_2 *\mathbb{Z}_2$, then $F$ is unique.
\end{lem}

\noindent {\bf Proof}  Note that both $\mathbb{Z}$ and $\mathbb{Z}_2 *\mathbb{Z}_2$ have no nontrivial finite normal subgroups.  If $F'$ was another finite normal subgroup of $H$, then $F'/(F' \cap F)$ would be a finite normal subgroup of $\mathbb{Z}$ or $\mathbb{Z}_2 *\mathbb{Z}_2$, thus must be trivial.  Hence $F' \subset F$.  In addition, if $H/F'$ were isomorphic to $\mathbb{Z}$ or $\mathbb{Z}_2 *\mathbb{Z}_2$, then we could make the same argument, reversing the roles of $F$ and $F'$, and get that $F = F'$.  \qed \\

We remark that the two possibilities in each of the third and fourth conditions of Theorem \ref{char_two-ended} above are distinguished by whether or not $H$ interchanges the two ``endpoints at infinity'' when acting on its own Cayley graph.  

They are also distinguished by whether or not the abelianization, $H_1(H)$, of $H$ is finite or infinite.  For if $H = F*_F$, then there is a surjection from $H_1(H)$ onto $\mathbb{Z}$, hence $H_1(H)$ is infinite.  If instead $H = A*_FB$, then $H_1(H)$ is a quotient of $A \times B$, hence is finite.

We now recall a fact from \cite{BH}:

\begin{thm} \label{finite_subgroups}
If $G$ is a hyperbolic group, then $G$ contains only finitely many conjugacy classes of finite subgroups.
\end{thm}

\noindent For a proof of this, see Theorem III.$\Gamma$.3.2 from \cite{BH}.

We will need one more standard lemma for the proof of Proposition \ref{elem}.  As we were unable to find a proof for it in the literature, we include one below:

\begin{lem} \label{stabs}
Let $G$ be a hyperbolic group, and let $N_1 \subset N_2 \subset N_3 \subset \ldots$ be an ascending union of two-ended subgroups of $G$.  Then the sequence must stabilize.
\end{lem}

\noindent {\bf Proof} We may assume that, for each $i$, $N_i$ is properly contained in $N_{i+1}$.  Then we must show that the sequence of $N_i$'s terminates.  Lemma \ref{F_unique} tells us that, for each $i$, $N_i$ contains a unique finite normal subgroup $F_i$ such that $N_i/F_i$ is isomorphic to $\mathbb{Z}$ or $\mathbb{Z}_2*\mathbb{Z}_2$.  In the former case, $H_1(N_i)$ is infinite, and in the latter case it is finite.  

Let $C_i$ denote the commutator subgroup of $N_i$, so that $H_1(N_i) = N_i/C_i$, and for all $i<j$, $C_i \subset C_j$.  Fix $i$ and $j$ with $i<j$, so that $N_i$ is a finite index subgroup of $N_j$.  It follows that if $C_i$ is a finite index subgroup of $N_i$, then $C_j$ is of finite index in $N_j$.
Thus, if $H_1(N_i)$ is finite, then $H_1(N_j)$ is finite, i.e. if $N_i/F_i \cong \mathbb{Z}_2*\mathbb{Z}_2$, then $N_{j}/F_{j} \cong \mathbb{Z}_2*\mathbb{Z}_2$.

Now, by perhaps removing finitely many subgroups from the beginning of the sequence of $N_i$'s, we may assume that either $N_i/F_i \cong \mathbb{Z}$ for all $i$, or that $N_i /F_i \cong \mathbb{Z}_2*\mathbb{Z}_2$ for all $i$.  We claim that $F_1 \subset F_2 \subset F_3 \subset \ldots$

In the case when each $N_i/F_i$ is isomorphic to $\mathbb{Z}$, this is immediate, since $F_i$ consists of all elements of $N_i$ which are of finite order.  So assume that each $N_i/F_i$ is isomorphic to $\mathbb{Z}_2*\mathbb{Z}_2$.  

Fix $i$, and let $A = N_i/F_i \cong \mathbb{Z}_2*\mathbb{Z}_2$, $A' = N_{i+1}/F_{i+1} \cong \mathbb{Z}_2*\mathbb{Z}_2$.  Then we have the following diagram, where the rows are exact, and $\iota$ is again inclusion:

\[ \begin{array} {ccccccccc}
1 & \to & F_i & \to & N_i & \to & A & \to & 1\\
&&&&\downarrow \iota \\
1 & \to & F_{i+1} & \to & N_{i+1} & \stackrel{p}{\to} & A' & \to &1\\
\end{array} \]
Note that the kernel of $p \circ \iota$ must be finite.  Since $N_i$ is virtually $\mathbb{Z}$, the image $p \circ \iota (N_i)$ is of finite index in $A'$, and thus must be isomorphic to $\mathbb{Z}$ or $\mathbb{Z}_2*\mathbb{Z}_2$.  Hence the kernel of $p \circ \iota$ must equal $F_i$, so $F_i $ is contained in the kernel of $p$, which is $F_{i+1}$.  It now follows that $F_1 \subset F_2 \subset \ldots$.

By Theorem \ref{finite_subgroups}, the sequence $F_1 \subset F_2 \subset \ldots $ eventually stabilizes.  Thus, again by omitting finitely many of the first subgroups $N_1, N_2, \ldots , N_k$, we can assume that the sequence $N_1 \subset N_2 \subset \ldots$ is such that $F_1 = F_2 = F_3= \ldots$  Denote this finite group by $F$.

We diverge now to introduce the translation function $\tau$.  Any hyperbolic group $\Gamma$ has a well-defined translation function $\tau : \Gamma \to \mathbb{R}_{\geq 0}$ which, for any $\gamma \in \Gamma$, is defined by
$$\tau (\gamma) = \lim_{n \to \infty} \frac{1}{n}d(1, \gamma^n),$$
and which is nonzero for any $\gamma$ of infinite order.
It follows that, for any $\gamma$, $\tau (\gamma^m) = |m|\cdot \tau (\gamma )$, and that $\tau$ is constant on conjugacy classes.  (See  \cite{BH}, Section III.$\Gamma$.3.)
Proposition III.$\Gamma$.3.15 of \cite{BH} says that for any $R \in \mathbb{R}$, there exist only finitely many conjugacy classes $[\gamma ]$ of elements of $\Gamma$ such that $\tau(\gamma)<R$.  Thus we may define 
$$\tau_i = \min_{g \in N_i, o(g) = \infty} \{ \tau (g)\},$$
and let $g_i \in N_i$ be an element which realizes this minimum.

For each $i$, $N_i \subset N_{i+1}$, and hence $\tau_i \geq \tau_{i+1}$.  Thus to show that the sequence of $N_i$'s terminates, it suffices to show that $\tau_i > \tau_{i+1}$.

Assume first that each $N_i / F $ is isomorphic to $\mathbb{Z}$.  The inclusion $N_i/F \subset N_{i+1} /F$ is proper, hence there is some $g \in N_{i+1}$ such that, for some $k>1$, $g^kF = g_iF$.  
Thus for each integer $a$ there is some $f \in F$ such that $(g^k)^a = g_i^af$.  As $F$ is finite, there must exist distinct $a, b \in \mathbb{Z}_{>0}$ and some $f \in F$ such that $(g^k)^a = g_i^{a}f$ and $(g^k)^b= g_i^{b}f$, so $(g^k)^{b-a} = g_i^{b-a}$.
Thus $k \cdot |b-a| \cdot \tau (g) = \tau (g^{k \cdot (b-a)}) = \tau (g_i^{b-a}) = |b-a|\cdot \tau (g_i)$, and therefore $\tau(g) = \frac{1}{k}\tau_i$, and $\tau_i > \tau_{i+1}$.

Assume instead each $N_j/F$ is isomorphic to $\mathbb{Z}_2*\mathbb{Z}_2$.  Let $M_j$ denote the subgroup of $N_j$ that contains $F$ and is such that $M_j/F \cong \mathbb{Z}$ is of index two in $N_j/F$.  Note that $M_j$ contains all of the elements of $N_j$ that are of infinite order, so $g_j \in M_j$.  It also follows from this that $M_i = N_i \cap M_{i+1}$, and so $M_1 \subset M_2 \subset M_3 \subset \ldots$  Thus we may apply the analysis of the previous case, so we always have that $\tau_i > \tau_{i+1}$.
\qed \\

We can now give the
\\

\noindent {\bf Proof of Proposition \ref{elem}}  Part 1. of Definition \ref{def_of_elem} is immediate.

For 2., let $H$ be a two-ended subgroup in $\mathscr{C}$.  We will show that the commensurizer of $H$ in $G$ is the unique maximal two-ended subgroup which contains $H$.   
We note that if $H'$ is any two-ended subgroup of $G$ which contains $H$, then $H$ is of finite index in $H'$, and so $H'$ is contained in the commensurizer of $H$.
Now property 2. will follow once we have shown that Comm$_G(H)$ is two-ended, i.e. in $\mathscr{C}$.

Recall the definition of $\tau : \Gamma \to \mathbb{R}_{\geq 0}$ from the proof of the last lemma.  If $g \in $ Comm$_G(H)$, then for some nonzero $i$ and $j$, $g h^i g^{-1} = h^j$, where $h$ generates an infinite cyclic subgroup of $H$.  Since $\tau (\gamma^m) = |m|\cdot \tau (\gamma )$ for each $\gamma$, and $\tau$ is constant on conjugacy classes, we must have that $|i| = |j|$.  Hence $g$ is in the normalizer $N_G\langle h^i\rangle $, and so
$${\rm Comm}_G (H ) = \cup_{i=1}^\infty N_G\langle h^i\rangle  .$$
Note that, for all $i$, the centralizer $Z_G\langle h^i \rangle $ is a subgroup of index at most two in $N_G\langle h^i \rangle $.  Corollary III.$\Gamma$.3.10(2) of [BH] says that each $Z_G\langle h^i \rangle$ is virtually $\mathbb{Z}$, hence $N_G\langle h^i \rangle $ is as well.  Using this, we will now see that, for some $m_0$, 
$\cup_{i=1}^\infty N_G\langle h^i\rangle  = N_G\langle h^{m_0}\rangle $.

Note that, for any $i_0$ which divides $i$, $N_G\langle h^{i_0}\rangle \subset N_G\langle h^i\rangle$.  Thus 
$$\cup_{i=1}^\infty N_G\langle h^i\rangle = \cup_{i=1}^\infty N_G\langle h^{i!}\rangle ,$$
and the latter is an ascending union.
By Lemma \ref{stabs}, the terms must eventually stabilize, hence for some $m_0$, Comm$_G(H)$ equals $N_G\langle h^{m_0} \rangle$, and thus is two-ended.

Property 3. follows from the implication from Theorem \ref{finite_subgroups} that the finite subgroups of $G$ are of bounded order.
As for 4., if $C \in \mathscr{C}$ acts on a tree without fixing a point, then $C$ is 
virtually $\mathbb{Z}$, hence has an axis, so preserves two points in the boundary of the 
tree.

Property 5. follows from knowing, from the proof of 2. above, that any maximal, infinite $C \in \mathscr{C}$ is the commensurizer of any two-ended subgroup of $G$ contained in $C$.  In particular, $C = $ Comm$_G(C)$, and since the normalizer of a subgroup is contained in its commensurizer, 5. follows.
\qed
\\


We shall now define the notion of complexity used in \cite{DP}, and then carefully state their result.  Let $G$ be a finitely presented group, let $\mathscr{C}$ be a family of elementary subgroups of $G$, and note that $G$ is the 
fundamental group of a two-dimensional simplicial orbihedron $\Pi$ for which vertex 
stabilizers are in $\mathscr{C}$.  For any such $\Pi$, we define $T(\Pi)$ to be the 
number of faces of $\Pi$, and $b_1 (\Pi)$ to be the first Betti number of the space.  
Then we define the complexity of $\Pi$ to be 
$$c(\Pi) = (T(\Pi ), b_1(\Pi)).$$
The complexity of $G$ with respect to $\mathscr{C}$ is defined to be
$$c(G, \mathscr{C}) = c(G) = \min c(\Pi),$$
where the minimum is taken over all $\Pi$ with vertex groups elements of $\mathscr{C}$ 
and $G = \pi_1^{orb} (\Pi)$.  Lexicographical ordering is used.

Proposition 3.4 of \cite{DP} shows that if $c(G) = (0,0)$, then $G$ must be the 
fundamental group of a tree of groups (possibly just a vertex), with finite edge stabilizers, and vertex stabilizers in $\mathscr{C}$.  We are taking $G$ to be finitely presented, so we note that any such tree will be finite.
A group  is said to have a {\it dihedral} action on a tree if the group acts on the tree, has an axis, and elements of the group interchange the endpoints of the axis.
\cite{DP} proves

\begin{thm}[\cite{DP}] \label{thm_3.6}
Let $G$ be a finitely presented group, with $\mathscr{C}$ a family of elementary subgroups of $G$.  Suppose $G$ has a proper decomposition over $\mathscr{C}$, with $\tau$ the associated Bass-Serre tree, and suppose further that no $C \in \mathscr{C}$ admits a dihedral action on $\tau$.  Then there is a decomposition of $G$ over $\mathscr{C}$ with associated tree  $\tau'$ such that there is a natural $G$-equivariant map $\tau' \to \tau$ which does not collapse edges to vertices, 
and, if $\{ G_v\}$ denotes the vertex groups of $G \backslash \tau'$, then $\sum T(G_v) \leq T(G)$, and $\sup_v c(G_v, \mathscr{C} \cap G_v) < c(G, \mathscr{C})$.
\end{thm}

With $\mathscr{C}$ defined to be the finite and two-ended subgroups of a group $G$, as long as $G$ has no 2-torsion, it follows that the action of any $C \in \mathscr{C}$ on any $G$-tree $\tau$ is not dihedral.

To make use of this theorem, we will need several lemmas.  In the arguments that follow, we shall denote vertices and edges with lower case letters, and their stabilizers with the capitalizations of those letters.


\begin{lem}\label{pushforward}
Let $G$ be a finitely generated group, and $\mathscr{C}$ a family of subgroups of $G$ 
which is closed under conjugation and subgroups.  Suppose that $\phi : \tau' \to \tau$ is 
a $G$-map between $G$-trees with the action of $G$ on $\tau$ minimal, and all edge 
stabilizers in $\mathscr{C}$.  Moreover, suppose $\phi$ is such that, for each edge $e$ of $\tau'$, stab$(e)$ is contained in stab($\phi 
(e)$) with finite index.  Let $\Gamma' = G \backslash \tau'$, and $\Gamma = G \backslash 
\tau$, and suppose that the edge groups of $\Gamma$ are all finitely generated.
Then if $\Gamma'$ admits a proper refinement over $\mathscr{C}$, so does $\Gamma$, and 
the additional edge groups in the refinements are the same.
\end{lem}

From this, we immediately have

\begin{cor}\label{wmaxtarget_implies_wmaxsource}
If $G$, $\mathscr{C}$, $\Gamma$ and $\Gamma'$ are as in Theorem \ref{pushforward}, and $\Gamma$ is a maximal proper decomposition of $G$, then $\Gamma'$ must be maximal as well.
\end{cor}

\noindent {\bf Proof of Lemma \ref{pushforward}}  Recall that a Stallings' fold on a $G$-tree is a $G$-map which takes two adjacent edges $e$ and $f$, meeting at vertex $v$, with $e$ also incident to vertex $x$ and $f$ to vertex $y$, and identifies $e$ to $f$, such that $x$ is identified to $y$.  Different types of Stallings' folds, as described in \cite{BF}, correspond to whether and with what orientation $e$ is in the $G$-orbit of $f$, and whether $x$ and/or $y$ are in the $G$-orbit of $v$.  The group $G$ and the edge stabilizers of $\tau$ are all finitely generated, $\phi$ is a simplicial map between trees, and $\phi$ collapses no edge of $\tau'$ to a vertex, thus \cite{BF} tells us that $\phi$ must be a composition $\phi_n \circ \phi_{n-1} \circ \ldots \circ \phi_1$ of Stallings' folds $\{ \phi_i\}$.

We will show that if $\phi$ is a Stallings' fold, then a proper splitting of a vertex group of $\tau' $, which is compatible with $\Gamma'$,  induces a proper splitting of the image of the vertex group which is compatible with $\Gamma$, over the same edge group.  From this, it will follow that a proper refinement of $\Gamma'$ induces a proper refinement of $\Gamma$.

So assume that $\phi : \tau' \to \tau$ is a Stallings fold.  We use our notation from above, so that $\phi$ identifies $e$ to $f$ and $x$ to $y$, where $e$ and $f$ meet at the vertex $v \in \tau'$, and similarly, identifies $g\cdot e$ to $g \cdot f$ for each $g$ in $G$.
Let vertex $w \in \tau'$ be such that $W$, the stabilizer of $w$, admits a proper splitting over some $C \in \mathscr{C}$, which is compatible with $\Gamma'$.  Thus there exists a tree $\overline{\tau'}$ and a $G$-equivariant map $\zeta' : \overline{\tau'} \to \tau'$ which merely collapses each edge in the orbit of $c$ to a vertex in the orbit of $w$.  We would like to find a tree $\overline{\tau}$ such that there is a similar collapsing map $\zeta :\overline{\tau} \to \tau$, a fold $\overline{\phi}$ taking $\overline{\tau'}$ to $\overline{\tau}$, and such that the following diagram commutes:

\[ 
\begin{array}{ccc}

\overline{\tau'} & \stackrel{\overline{\phi}}{\rightarrow} & \overline{\tau}\\
\zeta ' \downarrow && \downarrow \zeta \\
\tau' & \stackrel{\phi}{\rightarrow} &\tau\\

\end{array} 
\] 

For our first case, assume that $w$ is not in the $G$-orbit of $v$, nor of $x$ nor $y$.  Then we may define $\overline{\phi}$ to identify $\zeta'^{-1}(e)$ to $\zeta'^{-1}(f)$, and $\zeta'^{-1}(g\cdot e)$ to $\zeta'^{-1}(g\cdot f)$, for each $g$ in $G$.  The edge $c$, as well as the edges in the $G$-orbit of $c$, are untouched by such a fold, so the above diagram must commute.  Also because no edge gets identified to $c$ or any of its translates, and because the refinement $\overline{\Gamma'}$ of $\Gamma'$ is proper, it follows that $\overline{\phi}$ induces a refinement of $\Gamma$ which is proper.  

Next, assume that $w$ is in the $G$-orbit of $x$, and not of $v$.  ($w$ may be in the orbit of $y$.)  Then we may again define $\overline{\phi}$ directly, taking that it identifies $\zeta'^{-1}(e)$ to $\zeta'^{-1}(f)$, and similarly for the $G$-orbits of $e$ and $f$.  Define the map $\overline{\phi}_*$ to take the stabilizer of any vertex or edge $z$ in $\tau'$ to the stabilizer of $\overline{\phi}(z)$, and let $a$ and $b$ be the vertices of $c$.  Recall that $A = $stab$(a)$, and so on.  Then in this case, $\overline{\phi}_*(C)=C$, while $A \subseteq \overline{\phi}_*(A)$ and $B \subseteq \overline{\phi}_*(B)$.  It follows again that $\overline{\Gamma}$ is a proper refinement of $\Gamma$ because $\overline{\Gamma'}$ is  a proper refinement of $\Gamma'$.  To see this, we note that 
if $C \hookrightarrow A$ and $C\hookrightarrow B$ are not isomorphisms, then neither are the new injections in $\overline{\tau}$.  If instead $g \in V$ takes $a$ to $b$, then $g$ will take $\overline{\phi}(a)$ to $\overline{\phi}(b)$.  Thus we have that $\overline{\Gamma}$ is a proper refinement of $\Gamma$.

It remains to consider the case in which $w$ is in the $G$-orbit of $v$.  Without loss of generality, we assume that $w=v$.  By abuse of notation, we will denote $\zeta'^{-1}(e)$ by $e$, and $\zeta'^{-1}(f)$ by  $f$.  Suppose that $e$ and $f$ are adjacent in $\overline{\tau'}$, so both contain either $a$ or $b$.  Here again, we may simply define $\overline{\phi}$ to identify $e$ to $f$, and extend equivariantly.  Then $\overline{\phi_*}$ takes $A$, $B$ and $C$ to themselves, and if there is some $g \in V$ which takes $a$ to $b$, then $g$ must also take $\overline{\phi}(a)$ to $\overline{\phi}(b)$.  Hence, this induced refinement $\overline{\Gamma}$ must be proper.  

So for our last case, assume that $w=v$, and that $e$ and $f$ are not adjacent in $\overline{\tau'}$.  
Without loss of generality, take that $e$ contains $a$ and $f$ contains $b$, i.e. $E\subseteq A$ and $F \subseteq B$.  Here, we will use our hypothesis that $E$ and $F$ are of finite index in $\phi_* (E) = \phi_*(F)$ to show that either $E \subseteq C$ or $F \subseteq C$.  If $E \subseteq C$, then we may alter $\overline{\tau'}$ by `sliding' $e$ so that it is incident to $b$ instead of $a$, and do the same with the $G$-orbit of $e$.  If $F \subseteq C$, then we can slide $f$ instead.  By doing this, we are able to create a proper refinement of $\Gamma'$ of the type discussed in the previous paragraph, and may refer now to that argument.

To show that this is possible, assume that neither $E$ nor $F$ is contained in $C$, and choose elements $g_E \in E-C$ and $g_F \in F-C$.  Then the subset of $\overline{\tau'}$ which is fixed pointwise by $g_E$ is a subtree of $\overline{\tau'}$ which is disjoint from the subtree of points fixed by $g_F$.  Thus $g_Eg_F$ acts by translation on an axis in $\overline{\tau'}$.  
Both $E$ and $F$ are contained in $\phi_*(E)$, hence so is $g_Eg_F$, but because $g_Eg_F$ has an axis, it is of infinite order, and no power $(g_Eg_F)^n$ is contained in $E$ or $F$, except when $n=0$.  This means that $E$ and $F$ must be of infinite index in $\phi_*(E)$, which is a contradiction.  Thus either $E \subseteq C$ or $F \subseteq C$ as desired.

We have seen now that if $\phi : \Gamma' \to \Gamma$ is a fold, and  if $\Gamma'$ admits a proper refinement by a splitting over a subgroup $C$, then $\Gamma$ must also admit a proper refinement by a splitting which is also over $C$.  It follows that if $\phi$ is a composition of folds, then a proper refinement of $\Gamma'$ pushes through each fold, giving a proper refinement of $\Gamma$, as desired.  \qed
\\


In preparation for our proof of Lemma \ref{vtx_gps}, we will now present an analysis of Stallings' folds of type ``A''.

Assume that we have a Stallings' fold between $G$-trees, $\psi : \sigma' \to \sigma$, which takes edges $e$ and $f$ of $\sigma'$ which share a vertex $v$, and identifies $e$ to $f$, such that the other vertex, call it $x$, of $e$ is identified to the other vertex, $y$, of $f$, and makes similar identifications to $g\cdot (e \cup f)$, for all $g \in G$.  Assume further that neither $x$ nor $y$ are in the $G$-orbit of $v$. Then  $\psi$ must in fact be a fold of one of three types which, following \cite{BF}, we will call types IA, IIA, and IIIA.
These types correspond to the following three cases:  when no $g \in G$ takes $x$ to $y$, when some $g \in G$ takes $x$ to $y$ and $e$ to $f$, and when some $g \in G$ takes $x$ to $y$, but does not take $e$ to $f$.  

Let $\pi$ denote the projection maps $\sigma \to \Sigma = G \backslash \sigma$, and $\sigma' \to \Sigma' = G\backslash \sigma'$, and let $\Psi: \Sigma' \to \Sigma$ be such that $\Psi \circ \pi = \pi \circ \psi$.  Our figures below indicate how, in each case, $\Psi$ will alter $\pi (e \cup f)$.  Since $\Psi$ cannot alter the underlying graph, or edge or vertex groups, of $\Sigma' - \pi (e \cup f)$, these must describe $\Psi$ completely.

\begin{figure}[!h] 
\begin{center}
 \includegraphics[width=3.56in]{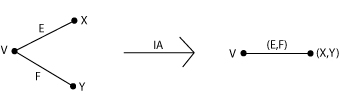} 
 \caption{A fold of type IA, with vertices and edges labeled with their associated groups}
 \label{IA}
\end{center}
\end{figure}

When no $g \in G$ takes $x$ to $y$, we will say that the fold is of type IA.  In this case, $\pi (e \cup f)$ will change as indicated in Figure \ref{IA}.

\begin{figure}[!h] 
\begin{center}
   \includegraphics[width=3.53in]{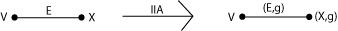} 
   \caption{A fold of type IIA}
\label{IIA}
\end{center}
\end{figure}

A fold of type IIA occurs when some $g \in G$ takes $x$ to $y$ and takes $e$ to $f$, in which case we have that $g \in V$, the stabilizer of $v$.  Here, the image under $\pi$ of the segment $e \cup f$ is a single edge, and folding changes only the labeling of $\Sigma'$.  See Figure \ref{IIA}.  

\begin{figure}[!h] 
\begin{center}
 \includegraphics[width=3.66in]{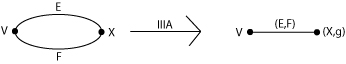} 
 \caption{A fold of type IIIA}
 \label{IIIA}
\end{center}
\end{figure}

Lastly, we have a fold of type IIIA when some $g \in G$ takes $x$ to $y$ and does not take $e$ to $f$.  Note that then $g$ translates along an axis containing $e$ and $f$.  In $\Sigma'$, we get what is shown in Figure \ref{IIIA}.

Next, we note the following fact, which we shall make use of with $n=2$:   if $\sigma$ is a $G$-tree, and $v_1, \ldots , v_n$ vertices of $\sigma$ with respective stabilizers $V_1, \ldots , V_n \subset G$, then the $\langle V_1, \ldots , V_n\rangle$-orbit of the smallest subtree containing $\{ v_1, \ldots , v_n\}$ is a $\langle V_1, \ldots , V_n\rangle$-tree.  This follows from 


\begin{lem}\label{orbit_makes_conn_subtree}
Let $G$ be a finitely generated group, with a $G$-tree $\sigma $ and associated decomposition $\Sigma$, identified with $G \backslash \sigma $.   Let $V_1, \ldots , V_n$ be stabilizers of vertices $v_1, \ldots , v_n$ of $\sigma $, and let $\sigma_0$ be the smallest subtree of $\sigma $ containing $\{ v_1, \ldots , v_n\}$.  Then the orbit of $\sigma_0$ under $\langle V_1, \ldots , V_n \rangle$ is connected, thus a subtree of $\sigma $.
\end{lem}

\noindent {\bf Proof:}  Fix any $w \in \langle V_1, \ldots , V_n\rangle$.  It will suffice to show that $w\cdot \sigma_0$ is connected to $\sigma_0$ in $\langle V_1, \ldots, V_n\rangle \cdot \sigma_0$.  

We can write $w = w_1w_2\cdots w_{m-1} w_m$, where each $w_i$ is contained in some $V_{j_i}$.  Then $w_m\cdot \sigma_0$ intersects $\sigma_0$ at the vertex stabilized by $V_{j_m}$, the subtree $w_{m-1}\cdot (w_m\cdot \sigma_0 \cup \sigma_0) = (w_{m-1}w_m\cdot \sigma_0) \cup (w_{m-1}\cdot \sigma_0)$ intersects $\sigma_0$ at the vertex stabilized by $V_{j_{m-1}}$, the subtree $w_{m-2} \cdot ((w_{m-1}w_m\cdot \sigma_0) \cup (w_{m-1}\cdot \sigma_0) \cup \sigma_0) = (w_{m-2}w_{m-1}w_m\cdot \sigma_0) \cup (w_{m-2}w_{m-1}\cdot \sigma_0) \cup (w_{m-2}\cdot \sigma_0)$ intersects $\sigma_0$ at the vertex stabilized by $V_{j_{m-2}}$,
and so on.  Continuing in this manner, it follows that the translates $w\cdot \sigma_0 = w_1w_2\cdots w_m\cdot \sigma_0$, $w_1w_2\cdots w_{m-1}\cdot \sigma_0$, $\ldots$, $w_1w_2\cdot \sigma_0$, $w_1\cdot \sigma_0$, $\sigma_0$ make a subtree, hence $w\cdot \sigma_0$ is connected to $\sigma_0$ in $\langle V_1, \ldots , V_n \rangle \cdot \sigma_0$.
  \qed
\\

We can now prove the following:


\begin{lem}\label{vtx_gps}
Let $\Gamma$ be a maximal proper decomposition of a finitely presented group $G$ over a family $\mathscr{C}$ which is closed under conjugation and subgroups.  Let $\Gamma'$ be the decomposition from Theorem \ref{thm_3.6} of \cite{DP},
and assume that edge groups of $\Gamma'$ are of finite index in edge groups of $\Gamma$, as in Theorem \ref{pushforward}.  Then, for each vertex group $V$ of $\Gamma$,  either $V$ is a vertex group of $\Gamma'$, or $V\in \mathscr{C}$.
\end{lem}

\noindent {\bf Proof}  Let $\tau$ be the Bass-Serre tree associated to $\Gamma$, and $\tau'$ to $\Gamma'$.  We may subdivide the  edges of $\tau$ and $\tau'$ so that, for each edge of $\tau$ and $\tau'$, the vertices of that edge are in different $G$-orbits, yet still $\phi : \tau' \to \tau$ is a $G$-map.  
Again, we recall from \cite{BF} that $\phi$ is a composition of folds.  Our subdivision of the edges of $\tau$ and $\tau'$ ensures that $\phi$ is, in fact, a composition of folds of types IA, IIA, and IIIA.

Assume first that $\phi$ is a fold of type IA, IIA, or IIIA.  Then, using that $\Gamma$ is maximal, we will show that, for any vertex group $Z$ of $\Gamma$, either $Z$ is isomorphic by the given injection to one of its edge groups, or $Z$ is a vertex group of $\Gamma'$, hence has smaller complexity than $G$.  Thus for a composition of such folds, a vertex group of the target decomposition is either a vertex group of the source decomposition, or is in $\mathscr{C}$.

We employ our previous notation, so that $\phi$ is a fold which takes edge $e$ of $\tau'$ to edge $f$, and vertex $x$ to vertex $y$, with $e$ and $f$ sharing the additional vertex $v$.  It is immediate that, for all vertices $z'$ of $\tau'$, stab($z'$) $ = $ stab($\phi (z')$) if $z'$ is not in the $G$-orbit of $x$ or $y$.  Hence it suffices to show the above statement for $Z=$ stab($\phi (x)$).  

\begin{figure}[!h] 
\begin{center}
 \includegraphics[width=1.31in]{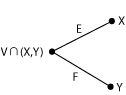}
  \caption{Decomposition of $Z$ in the case when $\phi$ is a fold of type IA.}
  \label{dIA}
\end{center}
\end{figure}

Consider the case in which $\phi$ is a fold of type IA.  Recall that $\phi (x) = \phi(y)$ has stabilizer $Z = (X,Y)$, and consider the action of $Z$ on $\tau'$.  Lemma \ref{orbit_makes_conn_subtree} implies that this gives the decomposition of $Z$ that is pictured in Figure \ref{dIA}.
Thus $(X,Y) = X *_E(V \cap (X,Y))*_FY$.  
If this decomposition gives a proper splitting of $Z$ which is compatible with $\Gamma$, i.e. the edge stabilizer of any edge adjacent to $\phi(x)$ is contained in a vertex group of the splitting, then this splitting would induce a proper refinement of $\Gamma$.  This would be a contradiction, however, because $\Gamma$ is assumed to be maximal.

We claim first that the decomposition is compatible with $\Gamma$, hence either splitting from the decomposition is compatible with $\Gamma$.  This follows because $E$ and $F$ are contained in $V$, so the stabilizer $(E,F)$ of $\phi(e)$ is contained in $V \cap (X,Y)$, and any other edge incident to $\phi (x)$ is untouched by the fold, hence has stabilizer either contained in $X$ or contained in $Y$.  

Therefore this decomposition of $Z$ must not give a proper splitting.  $[X*_E(V \cap (X,Y))]*_FY$ not being proper implies that either $Z=Y$ or $Y=F$.  If $Z=Y$, then $Z$ is a vertex group of $\tau'$.  Otherwise, $Y=F$.  But also $X*_E[(V \cap (X,Y)) *_FY]$ is not a proper splitting, so either $Z=X$ or $X=E$.  If $Z=X$, then, as before, $Z$ is a vertex group of $\tau'$.  Otherwise, $Y=F$ and $X=E$, so $Z = (X,Y) = (E,F)$, and hence $Z$ is an edge group of $\tau$.  Thus if $\phi$ is a fold of type $IA$, then either $Z$ is isomorphic to a vertex group of $\Gamma'$, or an edge group of $\Gamma$.

\begin{figure}[!h] 
\begin{center}
 \includegraphics[width=1.41in]{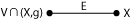} 
  \caption{Decomposition of $Z$ in the case when $\phi$ is a fold of type IIA.}
  \label{dIIA}
\end{center}
\end{figure}

Consider next the case in which $\phi$ is a fold of type IIA.  There is some $g \in G$ taking $e$ to $f$, and fixing $v$, and $\phi(x)$ is stabilized by $(X,g)$.  The action of this subgroup on $\tau'$ gives the splitting of $Z=(X,g)$ that is in Figure \ref{dIIA},
and hence $(X,g) = (V \cap (X,g))*_EX$.  It is clear that $(E,g) \subset (V \cap (X,g))$, so if we show that any other edge group of $\Gamma$ contained in $(X,g)$ is contained in one of the new vertex groups, then the compatibility of this splitting of $(X,g)$ with $\Gamma$ will follow.  But as above, since the fold only affects the edge group labeled $(E,g)$, then any other edge group incident to the vertex labeled $(X,g)$ must have been contained in $X$.

We note that since $g \in V \cap (E,g)$, but $g \notin E$, this splitting induces a proper refinement of $\Gamma$ unless $X=E$, in which case $(X,g) = (E,g)$.  Thus if $\phi$ is of type IIA, $Z = (X,g)$ must be an edge group of $\Gamma$.

\begin{figure}[!h] 
\begin{center}
 \includegraphics[width=1.41in]{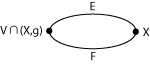} 
  \caption{Decomposition of $Z$ in the case when $\phi$ is a fold of type IIIA.}
  \label{dIIIA}
\end{center}
\end{figure}

If $\phi$ is a fold of type IIIA, then there is some $g \in G$ taking $x$ to $y$, but not taking $e$ to $f$.  Recall that $Z=$ stab($\phi(x)$) is $(X,g)$, and consider the action of $(X,g)$ on $\tau'$.  The quotient by this action contains the decomposition of $(X,g)$ given in Figure \ref{dIIIA},
thus $(X,g) = ((V \cap (X,g))*_FX)*_E$, where this HNN extension is by $g$.  A refinement by an HNN extension must always be proper, so it remains to show that this splitting induces a refinement of $\Gamma$, i.e. is compatible with the other splittings of $\Gamma$.  To do this, we must show that the stabilizer of any edge incident to $\phi (x)$ is contained in $((V \cap (X,g))*_FX)$.  The argument for this is similar to the above:  except for $\phi (e)$, any edge $d$ incident to $\phi (x)$ is again untouched by the fold, hence has stabilizer equal to the stabilizer of $\phi^{-1}(d)$, which is contained in $X$, as $\phi^{-1}(d)$ is incident to $x$.  $X \subset ((V \cap (X,g))*_FX)$, so our splitting is compatible with the splitting over $D$.  

Now recall that $\phi (e)$ is stabilized by $(E,F)$.  
But both $E$ and $F$ stabilize $v$, hence are in $V$.
Also, $E$ and $F$, when conjugated by $g$, stabilize $x$, hence $(E,F)$ is in $(X,g)$.  
Thus stab$(\phi (e)) = (E,F) \subset (V \cap (X,g)) \subset ((V\cap (X,g))*_FX)$, so the given splitting of $(X,g)$ is compatible with the other splittings of $\Gamma$.  But this means that there is a proper refinement of $\Gamma$, a contradiction.  Hence $\phi$ cannot be a fold of type IIIA.

We now address the situation in which $\phi = \phi_n \circ \phi_{n-1} \circ \ldots \circ \phi_1$, where each $\phi_i$ is a fold of type IA, IIA, or IIIA.  Let $\Gamma_i$ denote the decomposition $G \backslash \phi_i \circ \phi_{i-1} \circ \ldots \circ \phi_1 (\tau')$.  Theorem \ref{pushforward}, and the fact that $\Gamma$ is maximal, imply that the decompositions $\Gamma'$, $\Gamma_1$, $\Gamma_2$, $\ldots$, $\Gamma_{n-1}$ are all maximal.  Thus, for each $i$, the vertex groups of $\Gamma_i$ are edge groups of $\Gamma_i$, or are vertex groups of $\Gamma_{i-1}$.  It follows that any vertex group of $\Gamma$ is isomorphic to  either a vertex group of $\Gamma'$, or an edge group of some $\Gamma_i$, thus is in $\mathscr{C}$.  (Note that our early subdivision of edges of $\Gamma'$ only adds edge groups to the collection of vertex groups of $\Gamma'$, hence does not affect this result.)
\qed
\\


A result from \cite{SS} gives us the following theorem about the existence of maximal decompositions over two-ended subgroups:

\begin{thm}\label{max_decomps_exist}
Let $G$ be a one-ended, finitely presented group, and let $\Gamma$ be a proper decomposition of $G$ over two-ended subgroups.  Then $\Gamma$ admits a refinement $\Sigma$ which is a maximal proper decomposition of $G$ over two-ended subgroups.
\end{thm}

\noindent {\bf Proof:}  
 Let $\tau$ be the $G$-tree corresponding to $\Gamma$.  For any vertex $v$ of valence two of $\Gamma$ which is not the vertex of a circuit and has incident edges $e$ and $f$ such that $E=V=F$ by the given injections, collapse either $e$ or $f$.  Continue this process until no such vertices remain, and denote the resulting decomposition by $\overline{\Gamma}$.  
 We have now removed enough redundancy from $\Gamma$ to be able to apply
 Theorem 7.11 of \cite{SS}, with corrected statement in \cite{SSerrata}, giving us that $\overline{\Gamma}$ has a maximal refinement $\overline{\Sigma}$.
 
 We claim now that $\overline{\Sigma}$ induces a maximal refinement $\Sigma$ of $\Gamma$, i.e. that we may put the collapsed edges back into $\overline{\Sigma}$ corresponding to their location in $\Gamma$.  This can be done by merely subdividing each edge of $\overline{\Sigma}$ which corresponds to an edge $e$ (respectively $f$) of $\Gamma$ when, as in our notation above, the edge $f$ (respectively $e$) was collapsed to a point.
 \qed
\\


We can now prove Swarup's conjecture for hyperbolic groups with no 2-torsion:

\begin{thm} \label{VPC1acc_for_hyp} 
Let $G$ be a hyperbolic group with no 2-torsion.  Decompose $G$ maximally over finite subgroups, and 
then take the resulting vertex groups, and decompose those maximally over two-ended 
subgroups.  Now repeat this process on the new vertex groups and so on.  Then this 
process must eventually terminate, with subgroups of $G$ which are unsplittable over 
finite and two-ended subgroups.
\end{thm}

\noindent{\bf Proof}  We first will note that the above process must terminate for any finitely generated group $H$ such that $c(H) = (0,0)$, with respect to the family of finite and two-ended subgroups of $H$.  Recall that in this case, $H$ is the fundamental group of a tree of groups with finite edge stabilizers, and finite or two-ended vertex stabilizers.  If the tree consists of just one vertex, then $H$ is finite or two-ended.  When $H$ is finite, then it is unsplittable over all subgroups and hence the process terminates.  If $H$ is two-ended, then $H$ admits one nontrivial decomposition, which is over a finite subgroup and has finite vertex stabilizers, thus the above process must also terminate.

More generally, let $H$ be the fundamental group of a tree of groups as described above.  
The only vertex groups of the tree which admit any splittings are the two-ended groups.  As noted above, each splits over a finite subgroup, and the resulting vertex groups are finite, hence unsplittable.  Any collection of splittings of $H$ over finite subgroups are compatible, hence we may combine any splittings of vertex groups of the tree with the splittings of $H$ determined by the edges of the tree to get a decomposition of $H$ over finite subgroups with vertex groups which are completely unsplittable.  It follows that the process terminates for any $H$ such that $c(H)=(0,0)$.

Now we let $G$ be any hyperbolic group.  Then $G$ must be finitely presented, thus, by \cite{Dun}, it has a maximal decomposition over finite subgroups.  Choose such a decomposition, and let $\tau$ be the associated tree.  Let $\mathscr{C}$ be the family of all finite and two-ended subgroups of $G$, and let $\tau'$ be the $G$-tree resulting from an application of Theorem \ref{thm_3.6} of \cite{DP}.  Note that since the map $\tau' \to \tau$ collapses no edges to vertices, stabilizers of edges of $\tau'$ are subgroups of stabilizers of edges of $\tau$, hence the decomposition of $G$ associated to $\tau'$ is over finite subgroups of $G$.

Thus Lemma \ref{vtx_gps}, applied taking the family of elementary subgroups to be the collection of finite subgroups of $G$, implies that any vertex stabilizer $V_1$ of $\tau$ either is finite or is a vertex stabilizer of $\tau'$, hence is of smaller complexity (with respect to the family of finite and two-ended subgroups of $V_1$) than $G$.  Certainly the process described above must terminate for finite groups, so we may assume that $V_1$ is not finite.

Let $\mathscr{C}_1$ be the collection of finite and two-ended subgroups of $V_1$, i.e. $\mathscr{C}_1 = \mathscr{C} \cap V_1$.  Proposition \ref{elem} shows that $\mathscr{C}$ is elementary in $G$, and from this, it follows that $\mathscr{C}_1$ is elementary in $V_1$:  property 1. of Definition \ref{def_of_elem} is again immediate, and 4. follows from it holding in $\mathscr{C}$.  Finite subgroups of $G$ are of bounded order, hence finite subgroups of $H$ are also of bounded order, so 3. follows.  As for 2., we note that, for any two-ended $C \in \mathscr{C}_1$, any two-ended $C'$ containing $C$ must commensurize it, and $C \subset $ Comm$_{V_1}(C) \subset $ Comm$_G(C)$, thus Comm$_{V_1}(C)$ is two-ended, so again, Comm$_{V_1}(C)$ is the unique maximal element of $\mathscr{C}_1$ containing $C$.
Property 5. follows by the same argument as was given in Proposition \ref{elem}.  

Note that $V_1$ must have one end, so by Lemma \ref{max_decomps_exist}, $V_1$ has a maximal decomposition over two-ended subgroups. 
 Let $\tau_1$ be the corresponding $V_1$-tree, and $\tau_1'$ the tree from Theorem \ref{thm_3.6} of \cite{DP}.  
Since $V_1$ has one end, the edge groups of $\tau'$ are two-ended, and 
thus any edge group of $\tau'$ is of finite index in the image edge group from the map 
$\tau' \to \tau$.  Therefore, Lemma \ref{vtx_gps} gives us that if $V_2$ is a resulting 
vertex group, then $V_2$ is in $\mathscr{C}_1$ or has smaller complexity than $V_1$, with 
respect to the family $\mathscr{C}_2 = \mathscr{C}_1 \cap V_2$ of the finite and 
two-ended subgroups of $V_2$.  We note that $\mathscr{C}_2$ is also elementary in $V_2$.

If $V_2$ is in $\mathscr{C}_1$, then $V_2$ could admit one nontrivial maximal 
decomposition, which would be over a finite subgroup and would have finite vertex stabilizers.  Otherwise, we can repeat the arguments above, decomposing $V_2$ maximally over finite subgroups, decomposing the resulting vertex groups maximally over 2-ended subgroups, etc.
Complexity of the resulting groups continues to decrease, so we must eventually reach a collection of subgroups of $G$ which are unsplittable over any finite or two-ended subgroups, as desired.  \qed
\\


We will now use this result to get the hierarchy theorem for 3-manifolds stated earlier.  First, we recall that a 
surface $N$ in a 3-manifold $M$ is said to be essential if $N$ is properly embedded in 
$M$, 2-sided, $\pi_1$-injective into $M$, and is not properly homotopic into the boundary 
of $M$.

\begin{lem}
Let $M$ be a connected $3$-manifold with boundary, and let $\mathscr{A} = \{ A_i\}_{i \in 
I}$ be a finite collection of disjoint, non-parallel, essential
surfaces in $M$, such that $\{ \pi_1 (A_i)\}$ are contained in a family $\mathscr{C}$ of 
subgroups of $G = \pi_1(M)$ which is closed under subgroups and conjugation.  Suppose 
further that $\mathscr{A}$ is maximal with respect to collections of disjoint, 
non-parallel essential surfaces of $M$ with fundamental groups in $\mathscr{C}$.  Let 
$\Gamma$ be the decomposition of $G$ which is dual to $\mathscr{A}$.  Then $\Gamma$ is a 
maximal proper decomposition of $G$ over $\mathscr{C}$.
\end{lem}

\noindent {\bf Proof:}  Assume for the contrapositive that $\Gamma$ is not maximal.  Then 
there exists some vertex group $V$ of $\Gamma$ which admits a proper splitting over some 
$C \in \mathscr{C}$ which is compatible with $\Gamma$.  Let $L$ be the graph of groups 
for such a splitting of $V$, and let $p$ denote the midpoint of the edge of  $L$.
Let $N$ denote the union of the component of $M - \mathscr{A}$ which corresponds to $V$ 
with the surfaces $A_i$ which correspond to the edge groups incident to $V$.  
We can define a map from $N$ to $L$ which is an isomorphism on $\pi_1$, with each $A_i$ 
in $N$ mapped to a vertex of $L$, and such that the map is transverse to $p$.  Note that 
each component of the inverse image of $p$ is a properly embedded, 2-sided surface in 
$N$.  Furthermore, Stallings showed in \cite{Stl} that we can homotope this map on $N$ 
rel $\partial N$ to a new map $f$ such that the surfaces comprising $f^{-1}(p)$ are 
$\pi_1$-injective in $M$ (see also \cite{He}).

We may further assume that these components are not parallel to the boundary of $N$, 
because of the following.  Let $S$ denote a component of $f^{-1}(p)$ which is boundary 
parallel in $N$, and let $R$ be the region made up of $S$ and the component of $N-S$  
through which $S$ can be homotoped to $\partial N$, so $R$ is homeomorphic to $S \times 
I$.  Then we may homotope $f$ to take $R$ to $p$, and then to take a small neighborhood 
of $R$ past $p$, so that $p$ is not contained in $f(R)$.  We may then homotope $f$ to map 
the elements of $\mathscr{A} \cap R$ to the other vertex of $L$, so that still $p$ is not in 
$f(R)$, and still $f$ is an isomorphism on $\pi_1$.
Note that, because $L$ is the graph of groups of a proper splitting, and $f$ is 
surjective on $\pi_1$, this process will never make $f^{-1}(p)$ empty.

We have arranged that the components of $f^{-1}(p)$ are essential in $M$.  Because $f$ is 
$\pi_1$-injective, the fundamental group of each component of $f^{-1}(p)$ is conjugate to 
a subgroup of $C$ and so is in $\mathscr{C}$.  Since $f$ maps the $A_i$'s to vertices of 
$L$, the surfaces $f^{-1}(p)$ are disjoint from $\mathscr{A}$.  Also, as components of 
$f^{-1}(p)$ are not boundary parallel in $N$, they are not parallel to elements of 
$\mathscr{A}$.  Hence $\mathscr{A}$ is not maximal.  \qed \\

We note that each component of $f^{-1}(p)$ induces a refinement of $\Gamma$.  Suppose, in 
addition to the hypotheses on $M$ in the above lemma, that $M$ is irreducible.  Then we can 
homotope $f$ to remove any sphere components of $f^{-1}(p)$, so that any simply connected 
component of $f^{-1}(p)$ must be a compressing disk for $M$.  Thus, a maximal collection 
of compressing disks in an irreducible, connected 3-manifold $M$ induces a maximal proper 
decomposition of $G$ over $\{ 1\}$.

It also follows that, if $\mathscr{A}$ is a maximal collection of annuli in $M$, and $M$ 
is as in the above lemma, has incompressible boundary and is irreducible, then the graph 
of groups $\Gamma$ corresponding to $\mathscr{A}$ must be maximal over the family 
generated by all infinite cyclic subgroups of $\pi_1(M)$.

Recall that, if $M$ is orientable and irreducible and $\pi_1(M)=G$ is infinite, then 
$G$ has no torsion (\cite{He}).  Hence any essential surface in $M$ with finite 
fundamental group must be simply connected, and any essential surface with two-ended 
fundamental group must be an annulus.   These observations, together with Theorem 
\ref{VPC1acc_for_hyp}, imply the following theorem.  We note that $M$ as above has hyperbolic fundamental group if and only if $M$ is hyperbolic and has no torus boundary components.

\begin{thm} \label{3mans}
Let $M$ be an irreducible, orientable, compact 3-manifold with hyperbolic fundamental group.
The process of decomposing $M$ along any maximal, disjoint collection of compressing 
disks, then decomposing the resulting manifolds along maximal, disjoint collections of 
essential annuli, then the resulting manifolds along compressing disks, then again along 
essential annuli and so on, must eventually terminate with a collection of manifolds 
which have incompressible boundary and admit no essential annuli, or are 3-balls.
\end{thm}

{\bf Acknowledgment}
The author gratefully thanks her advisor, Peter Scott, for his guidance.

\bibliographystyle{alpha}
\bibliography{my_biblio}

\end{document}